\documentclass[12pt]{article}
\usepackage[cmtip,arrow]{xy}
\usepackage{pb-diagram,pb-xy}
\usepackage{amssymb} \usepackage{amsfonts}
\usepackage{amsmath}
\begin{document}
\title{Intersection forms of almost-flat 4-manifolds}
\author{A. Szczepa\'nski}
\date{\today}
\maketitle
\begin{abstract}
We calculate intersection forms of all 4-dimensional {\em almost-flat} manifolds. 
\end{abstract}
\newcommand{\F}{{\mathbb F}}
\newcommand{\Z}{{\mathbb Z}}
\newcommand{\Q}{{\mathbb Q}}
\newcommand{\R}{{\mathbb R}}
\newcommand{\C}{{\mathbb C}}
\newcommand{\h}{{\mathbb H}}
\newcommand{\N}{{\mathbb N}}
\section{Introduction}
Let $M$ be a smooth closed, oriented $4$-manifold
We define a symmetric and bilinear form (the intersection form)
$$Q_{M}:H^2(M,\Z)\times H^2(M,\Z)\to\Z,$$ 
\vskip 1mm
\noindent
as the evaluation of the cup product $a\cup b$ on the fundamental homology class of $[M]$, that is
$Q_{M}(a,b) = \langle a\cup b,[M]\rangle, a,b\in H^2(M,\Z).$
It remains a central question in 4-manifold topology which quadratic forms occur as the intersection form of an
orientable 4-manifold. In the topological category, there is the celebrated result of Freedman which asserts
that for each quadratic form $Q$ there exists an oriented simple-connected 4-manifold with $Q$ as its intersection form.
However in the smooth category Donaldson proved \cite{D} that among the definite quadratic forms only the diagonalilizable ones can be realized.
See \cite{D},\cite{DK},\cite{Scorpan}.
In this note we are interested in a classification of intersection forms of oriented $4$-dimensional {\em almost-flat} manifolds.
They are given by almost Bieberbach groups over 2 - and 3 - step nilpotent group.
We prove that if $M$ is not the torus $T^4$ then $Q_{M} = nH$ with $n = b_1(M) - 1,$ where $H$ is the hyperbolic form with
intersection matrix 
$\left[
\begin{smallmatrix}
0&1\\
1&0
\end{smallmatrix}\right]$
of rank 2 and $b_1(M)$ is the first Betti number of $M.$

An {\em almost-flat} manifold is a closed manifold $M$ such that for any $\epsilon > 0$ there exists
a Riemannian metric $g_{\epsilon}$ on $M$ with $|K_{\epsilon}|$diam($M,g_{\epsilon}$)$^2 <\epsilon$ where
$K_{\epsilon}$ is the sectional curvature and diam($M,g_{\epsilon}$) is the diameter of $M.$ 
When $K = 0$ $M$ is a flat manifold. The fundamental group of the almost-flat manifold is called
almost-Bieberbach group. Two almost-flat manifolds with isomorphic fundamental groups are affinely diffeomorphic,
see \cite[page 16]{dekimpe}.
From now on $M$ stands for an oriented $4$-dimensional {\em almost-flat} manifold. It is well known that 
the Euler characteristic of $M$ is zero (cf. \cite[p.134]{dekimpe}). 
Hence, from Poincare Duality we have the following relation involving the Betti numbers of $M$ 
\begin{equation}\label{betti}
\chi(M) = 2 - 2b_{1}(M) + b_{2}(M) = 0
\end{equation}
and we can assume that $b_{1}(M)\geq 1.$
Moreover we have
\newtheorem{lem}{Lemma}
\begin{lem}
For any oriented {\em almost-flat} manifold $M,$ the intersection form $Q_{M}$ is even.
\end{lem}
\vskip 1mm
\noindent
{\bf Proof:} From \cite[Lemma 1]{Bohr}, (see also the Wu's formula \cite[Theorem 11.14]{MS}) 
any closed, oriented $4$-dimensional spin-manifold has even intersection form.
Let us assume that $M$ has not a spin-structure. In \cite{LPS} and \cite{PS} all such manifolds
are classified. That means there is a list of their fundamental groups.
We claim that the first Betti number of any $M$ or equivalently
the rank of an abelianzation of $\pi_{1}(M)$ is equal to $1.$
For the proof we can use two methods. The first one uses presentations of $\pi_{1}(M)$
and computer system GAP, \cite{gap}. The second one uses properties of fundamental groups
and is as follows.
There are 7 flat manifolds without spin structure. All of them are presented in \cite{PS}.
Let us recall (see \cite{S}) that the fundamental group $\Gamma$ of a four dimensional flat manifold is the middle term in short exact
sequence
\begin{equation}\label{short}
0\to\Z^4\to\Gamma\stackrel{p}\to G\to 0,
\end{equation}
where $\Z^4$ is maximal abelian torsion free abelian group of rank 4 and $G$ is a finite group.
It is easy to see (c.f.\cite[page 51]{S}) that the first Betti number of the group $\Gamma$ is equal to the rank $(\Z^4)^{G},$
where the action of $G$ on $\Z^4$ is the following 
$$\forall g\in G, \forall z\in\Z^4 gz = \bar{g}z\bar{g}^{-1}.$$
Here $\bar{g}\in\Gamma$ is such that $p(\bar{g}) = g.$

In the {\em almost-flat} case we have a classification of all four dimensional oriented manifolds without spin structure \cite[p.12]{LPS}.
Let $E$ be a fundamental group of such manifold.
It is well known \cite{dekimpe} that there is the following generalization of the short exact sequence (\ref{short})
$$0\to N\to E\to G\to 0,$$ where $N$ is a nilpotent group and $G$ has a finite order.
In our situation, we can consider only a case where $N$ is 2 - or 3 - step nilpotent.
In the 2 - step nilpotent case, (see \cite[Theorem 6.4.10]{dekimpe}) the first Betti number of $E$ is equal to the first Betti number of some
3-dimensional crystallographic group 
$Q = E/\sqrt[N]{[N,N]},$ where
$$\sqrt[N]{[N,N]} = \{n\in N\mid n^k\in [N,N]\hskip 2mm \text{for some $k\geq 1$}\}.$$ 
The group $Q$ is called the underlying crystallographic group of $E.$  
In the 3 - step nilpotent case, also abelianization depends on the underlying crystallographic group. 
In order to obtain our results, we can apply
\cite[Remark 6.4.16]{dekimpe} and matrices on pages 225 - 230 of \cite{dekimpe}. 
For example, let us consider the almost Bieberbach group E with holonomy group $\Z_2$ (case 5 from
\cite[page 171]{dekimpe}) with the underlying crystallographic group $Q = C2.$ From definition $N$ is 2 - step nilpotent. 
Hence $C2$ denotes a 3-dimensional crystallographic group.
It is easy to see (cf.\cite[Table 3]{RS}) that $H_1(Q,\Z)\simeq \Z\oplus\Z_{2}^{2}$ and $b_{1}(E) =1.$
\footnote{In \cite[page 794]{GNA} a calculation of the abelianization of this group is not correct
but this does not effect the claim in the example or its proof.}
Hence from (\ref{betti}) $Q_{M} = 0.$ 
\vskip 2mm
\hskip 120mm
$\Box$

\vskip 5mm
\newtheorem{theo}{Theorem}
\begin{theo}
Let $M$ be any {\em almost-flat} oriented $4$-manifold different from the torus with intersection form $Q_{M}$.
Then 

$$Q_{M} = \left\{ \begin{array}{lll}
0& \mbox{for $b_{1}(M) = 1$}\\
H& \mbox{for $b_{1}(M) = 2$}\\
2H& \mbox{for $b_{1}(M) = 3$}\end{array}\right.$$.
\end{theo}
\vskip 1mm
\noindent 
{\bf Proof:} From Lemma 1 $Q_{M}$ is even and from \cite{dekimpe} $b_1(M)\leq 3.$  
Hence it follows from \cite[Theorem 3]{Kim} $Q_{M} = nH.$  
An application of the formula (\ref{betti}) gives us the equation $n = b_{1}(M) - 1.$
\vskip 2mm
\hskip 120mm
$\Box$

\vskip 4mm
\noindent\textbf{Acknowledgment}
\vskip2mm
\noindent
We would like to thank K. Dekimpe, R. Lutowski and M. Mroczkowski for some useful comments.

{Institute of Mathematics, University of Gda\'nsk}\\
{ul. Wita Stwosza 57,\\
80-952 Gda\'nsk,\\
Poland}\\
{E-mail: matas@ug.edu.pl}


\begin{thebibliography}{99}
\bibitem{Bohr} Ch. Bohr, On the signatures of even $4$-manifolds, Math. Proc. Cambridge, 132 (2002), no.3, 453 - 469
\bibitem{dekimpe} Dekimpe, K., \emph{Almost-Bieberbach Groups : Affine and Polynomial Structures},
Lecture Notes in Mathematics 1639, Springer (1996).
\bibitem{D} S. Donaldson, An application of gauge theory to four dimensional topology. Journal of Differential Geometry,
18, (1983), 279-315
\bibitem{DK} S. Donaldson, P. Kronheimer. \emph{The Geometry of Four-Manifolds}, Oxford University Press, Oxford, 1991
\bibitem{gap} The GAP Group, GAP -- Groups, Algorithms, and Programming, Version 4.4.12, 2008, (http://www.gap-system.org)
\bibitem{GNA} A. G\c{a}sior, N. Petrosyan, A, Szczepa\'nski, Spin structures on almost-flat manifolds. Algebr.
Geom. Topol., 16(2) (2016), 783-796.
\bibitem{Kim} Jin-Hong Kim, The $\frac{10}{8}$-conjecture and equivariant $e_{C}$-invariants, Math. Ann., 329, 31-47, (2004) 
\bibitem{LPS} R. Lutowski, N. Petrosyan, A. Szczepa\'nski, Classification of spin structures on $4$-dimensional almost-flat manifold, 
accepted to Mathematika, arXiv:1703.04972
\bibitem{MS} J. W. Milnor, J. D. Stasheff \emph{Charactersistic classes}, Princeton University Press, Princeton,
N. J. and University of Tokyo Press, Tokyo, 1974. Annals of Mathematics Studies, No. 76.
\bibitem{PS} B. Putrycz, A. Szczepa\'nski, Existence of spin structures on flat manifolds, Adv. Geometry 10 (2), (2010), 323-322
\bibitem{RS} J. Ratcliffe, S. Tschantz, Abelianization of space groups, Acta Crystallogr. 65 (1), (2009), 18-27
\bibitem{Scorpan} A. Scorpan \emph{The wild world of $4$-manifolds}, American Mathematical Society, 2005
\bibitem{S} A. Szczepa\'nski, \emph{Geometry of crystallographic groups}, World Scientifice, 2012
\end{thebibliography}
\end{document}